\documentclass[13pt]{article}
\usepackage[cp1251]{inputenc}
\usepackage[english]{babel}
\usepackage{amssymb,amsmath,amsthm,indentfirst}
\usepackage{multicol}
\usepackage{graphicx}

\begin{document}

\begin{center}
\LARGE{Non Controllability to Rest of the Two-Dimensional Distributed System Governed by
the Integrodifferential Equation}
\end{center}

\vskip 0.5cm
\begin{center}
{\large Igor Romanov}
\footnote{National Research University Higher School of Economics,\\
20 Myasnitskaya Ulitsa, Moscow 101000, Russia}
{\footnote{
E-mail: \it{ivromm1@gmail.com}}}
{\large Alexey Shamaev}
\footnote{Institute for Problems in Mechanics RAS,\\
101 Prosp. Vernadskogo, Block 1, Moscow 119526, Russia}
\footnote{Lomonosov Moscow State University,\\
GSP-1, Leninskie Gory, Moscow 119991, Russia}

\end{center}




\begin{abstract}
The paper deals with controllability problem for a distributed system
governed by the two-dimensional Gurtin-Pipkin equation. We consider a system
with compactly supported distributed control and show that if the
memory kernel is a twice continuously differentiable function, such that its
Laplace transformation has at least one non-zero root,
then the system cannot be driven to the equilibrium in a finite time.
\end{abstract}


\section{Introduction}

Integrodifferential equations with nonlocal terms of the convolution type often arise in applications such as mechanics of heterogeneous media, the theory of viscoelasticity, thermal physics, kinetic theory of gases and others.

For example, it was rigorously proved that in the case of heterogeneous two-phase medium,
consisting of viscous fluid and elastic inclusions, the effective equation is integro-differential,
and the corresponding convolution kernel is a finite of infinite sum of decreasing exponential
functions.

If the viscosity of the liquid is small (big), the effective equation does contain (does not contain) the third order terms corresponding to Kelvin-Voigh friction, see \cite{Sanchez}.

In the theory of viscoelasticity, it is a common practice to approximate the relaxation kernels
by the sum of exponents.

In thermal physics, laws of heat conduction with an integral memory are studied in many papers; note among them \cite{Gurtin}.

The presence of integral memory in the law of heat conduction might lead to the appearance
of a thermal front which moves at a finite speed. This makes an important difference with
the heat equation whose solution propagates at infinite speed.

In this paper, we give an outline of the results on the existence and uniqueness of solutions
to these systems and consider the problem of controllability.

\section{Statement of the Problem }
\label{intro}
In this article, we  consider the problem of non controllability of a system governed by
the integrodifferential equation
\begin{equation}
\label{1}
\theta_{t}(t,x,y)-\int\limits^t_0K(t-s)\Delta\theta(s,x,y)ds=u(t,x,y),
\end{equation}
$$
t>0,\ \ (x,y)\in\Omega.
$$
\begin{equation}
\label{2}
\theta|_{t=0}=\xi(x,y),
\end{equation}
\begin{equation}
\label{3}
\theta|_{\partial\Omega}=0,
\end{equation}
hereinafter
$\Omega\subset\mathbb{R}^2$ is a bounded domain, $K(t)$ is an arbitrary twice continuously differentiable
function such that $K(0)=\mu>0$, and
$u(t,x,y)$ is a control supported (in $x, y$) on $\Omega$.
The kernel $K(t)$ can be represented, for example, as a sum of decreasing exponential functions:
$$
K(t)=\sum\limits_{j=1}^{N}c_je^{-\gamma_jt},
$$
$c_j$, $\gamma_j$ are given positive constants.

For brevity, we  write $\theta(t)$ and $u(t)$ instead of $\theta(t,x,y)$ and $u(t,x,y)$, respectively.
This also means that $\theta(t)$, $u(t)$ are functions of $t$ with values in some suitable
space.

The goal of the control is to drive this mechanical system to rest in a finite time.
We say that the system (\ref{1})---(\ref{3})
is \emph{controllable to rest} if for every initial condition $\xi$ we can find a control $u$
with compact support (in $t$)
such that the corresponding solution $\theta(t,u)$ of the
problem (\ref{1})---(\ref{3}) has a compact support (in $t$).
Conversely, the system is \emph{uncontrollable to rest} if there is the initial condition $\xi$ such that
for every control $u$ ($u$ is in the suitable class of functions) the corresponding solution does not
have compact support (in $t$).

In this article, we prove that the system governed by two-dimensional Gurtin-Pipkin equation
is uncontrollable to rest if the distributed control is supported on the
subdomain which is properly contained in arbitrary bounded domain with a smooth boundary. This result is some
generalization of the analogous theorem in \cite{Ivanov} devoted to the similar
one-dimensional problem. The method used in the paper can also apply to
the case where the dimension of $\Omega$ is greater than 2.
It will be discussed in Section 6.

\section{Literature Review}

The presence of nonlocal terms of the convolution type in the equations and systems leads to a number of interesting qualitative effects that are not observed in the case of differential equations and systems of equations.  For instance, the systems of this type exhibit the properties of both parabolic and hyperbolic
equations.
In spectral problems for such equations and systems the spectrum is composed of two parts: real and complex.
The former one corresponds to the energy dissipation in the heat equation; the later corresponds to  vibrations.
Such equations can be solved by means of the method similar to the Fourier method.

In addition, systems of this type are usually uncontrollable to rest, if we apply boundary control
or control which is distributed on the part of the domain.
Here we recall the well-known work \cite{Butk} that deals with the equation of the vibration of a string. As was proved in this work, if we apply the control to the end of the string, then the system can be driven to rest.
The author used the so-called moment method.
The mentioned results were generalized to the multidimensional case in \cite{Lions}.

At the same time, if we use the control distributed on the whole domain then integral terms of the
convolution type "facilitate" the process of control. In this case, the control time is significantly reduced.
It should be noted that the spectral method proposed in \cite{Chernousko} can be successfully adapted  to the case of systems with nonlocal terms of convolution type, see \cite{Romanov}.

The uncontrollability mentioned above was justified in \cite{Ivanov}
for one-dimensional systems similar to (\ref{1}).
In most cases the property of controllability to rest is not observed.
For example, in \cite{Ivanov} it was proved that a solution to the heat equation with memory cannot be driven to rest in a finite time if some auxiliary function has roots. This result is valid both for
boundary and distributed control. Moreover, the case of distributed control can be reduced to the case
of boundary control.
In our paper, we obtain similar results for the case of two-dimensional domains.

We should also mention the work \cite{Imanuvilov}, where the boundary non controllability was justified
for the heat equation with memory.


Positive results on controllability of a one-dimensional
wave equation with memory were obtained in \cite{Romanov}.
It was shown that this equation can be driven to rest by applying a bounded
distributed control. In this case, the kernel of the integral term in the equation is the sum
of $N$ decreasing exponential functions.

Problems similar to (\ref{1})---(\ref{3}) for integrodifferential equations were widely studied in
the existing literature.
Equation (\ref{1}) was originally derived in \cite{Gurtin}.
The questions of solvability and asymptotic behavior of solutions for equations of this type  were investigated for example in \cite{Dafermos}, \cite{Desh}.
In \cite{Rivera} it was proved that the energy for
some dissipative system decays polynomially, when the memory kernel decays exponentially.

Problems of solvability of system (\ref{1})---(\ref{3}) were considered in \cite{Vlasov}.
It was proved that a solution belongs to some Sobolev space on the semi-axis (in $t$)
if the kernel $K(t)$ is the sum of exponential functions, each of them tends to
zero as $t\rightarrow+\infty$.

Interesting explicit formulas for the solution of (\ref{1})---(\ref{3})
were obtained in \cite{Rautian} under the assumption that the kernel $K(t)$ is also the sum of decreasing exponential functions.
It follows from these formulas that solutions tend to zero when $t\rightarrow+\infty$.
In all these works it is supposed that the kernels of integral terms in the studied equations are non increasing functions.

\section{Preliminaries}

Let $A:=\Delta$ be an operator acting on a space
$D(A)=H^2(\Omega)\cap H^1_0(\Omega)$ where $\Omega\subset R^2$ is a bounded domain
with boundary of class $C^2$. We consider now the control function $u(t)\in C(0,+\infty;L_2(\Omega))$
and the initial condition $\xi\in H^2(\Omega)\cap H^1_0(\Omega)$.

\textbf{Definition 1.} The function
$$
\theta(t)\in C^1(0,+\infty;L_2(\Omega))\cap C(0,+\infty;H^2(\Omega)\cap H^1_0(\Omega))
$$
is the solution of the problem (\ref{1})---(\ref{3}) if $\theta(t)$ satisfies the equation (\ref{1}) and
initial condition (\ref{2}).

We note that the boundary condition (\ref{3}) makes sense because for any $t\geqslant 0$ $\theta(t)$ is
a continuous (in $x,y$) function. There are several theorems of existence and uniqueness dedicated
to the problem (\ref{1})---(\ref{3}) (see \cite{Pandolfi}).

Let us denote $PW_+$ as the linear space of the Laplace transforms of
elements of $L_2(0,+\infty)$ with compact support
distributed on $[0,\infty)$. It is a well known fact that
$\varphi(\lambda)\in PW_+$ if and only if it is an entire function such that

1) there are real numbers $C$ and $T$ such that $|\varphi(\lambda)|\leq Ce^{T|\lambda|}$. Note that $C$ and $T$
depend on $\varphi(\lambda)$.

2) $\sup\limits_{x\geq 0}\:\int\limits_{-\infty}^{+\infty}|\varphi(x+iy)|^2dy<+\infty$.

\section{The main results}

Now we consider an auxiliary boundary-value problem
\begin{equation}
\label{4}
\theta_{t}(t,x,y)-\int\limits^t_0K(t-s)\Delta\theta(s,x,y)ds=0,
\end{equation}
$$
t>0,\ \ (x,y)\in\Omega_0=\{(x,y): x^2+y^2<R^2\},
$$
\begin{equation}
\label{5}
\theta|_{t=0}=\xi(x,y),
\end{equation}
\begin{equation}
\label{6}
\theta|_{\partial\Omega_0}=v(t,x,y),\ \ (x,y)\in\partial\Omega_0.
\end{equation}

In this problem for each $T>0$ $\xi\in H^2(\Omega_0)$,
$v\in C(0,T;H^{\frac{3}{2}}(\partial\Omega_0))$.

\textbf{Definition 2.}
$$
\theta(t)\in C^1(0,+\infty;L_2(\Omega_0))\cap C(0,+\infty;H^2(\Omega_0))
$$
is the solution of the problem (\ref{4})---(\ref{6}) if
$\theta(t)$ satisfies the equation (\ref{4}),
initial condition (\ref{5}) and the boundary condition (\ref{6}) (in the sense of the trace).

Suppose there is a solution to the problem (\ref{4})---(\ref{6}).
We multiply (in the sense of the inner product in $L_2(\Omega_0))$ both parts of
(\ref{4}) by the function $\varphi$ such that $\varphi\in H^2(\Omega_0)\cap H^1_0(\Omega_0)$.
Hereinafter $\nu$ is a normal vector to the domain boundary $\partial\Omega_0$. After that,
using Green's formula we replace the operator $A$ from $\theta(t)$ to $\varphi$.
\begin{equation}
\label{7}
\frac{d}{dt}\langle\theta(t),\varphi\rangle-
\int\limits^t_0K(t-s)\left(\langle\theta(s),\Delta\varphi\rangle-
\int\limits_{\partial\Omega_0}v(s)\frac{\partial\varphi}{\partial\nu} d\sigma\right)ds=0,
\end{equation}
where $\langle\cdot,\cdot\rangle$ is the inner product in $L_2(\Omega_0)$.

The orthonormalized system of eigenvectors of $A$ are the functions $\varphi_{nm}(x,y)$
which in polar coordinates $x=r\cos\alpha$, $y=r\sin\alpha$ have the form
$$
\tilde{\varphi}_{nm}(r,\alpha)=\frac{J_m\left(\mu^m_n\frac{r}{R}\right)e^{im\alpha}}
{\sqrt{\pi}RJ_m^{\prime}(\mu^m_n)},\quad m=0,1,2,...,\ \ n=1,2,...,
$$
where $J_m$ are Bessel functions, $\mu^m_n$ are positive roots of $J_m$.
It is a well-known fact that this system is
a basis for $L_2(\Omega_0)$. We substitute $\varphi=\varphi_{nm}$ in (\ref{7}). Then using
the notation $\theta_{nm}(t)=\langle\theta(t),\varphi_{nm}\rangle$ we obtain
\begin{equation}
\label{8}
\frac{d\theta_{nm}(t)}{dt}+\lambda^2_{nm}
\int\limits^t_0K(t-s)\theta_{nm}(s)
d\sigma ds=-\int\limits^t_0K(t-s)\left(
\int\limits_{\partial\Omega_0}v(s)\frac{\partial\varphi_{nm}}{\partial\nu}\right)d\sigma ds
\end{equation}
where $\lambda_{nm}$ are the corresponding eigenvalues. If we use polar coordinates we have
$$
\lambda^2_{nm}=\frac{\left(\mu^m_n\right)^2}{R^2}.
$$

Let us make the Laplace transformation for both parts of (\ref{8}) and express $\hat{\theta}_{nm}(\lambda)$:
\begin{equation}
\label{9}
\hat{\theta}_{nm}(\lambda)=
\frac{-\hat{K}(\lambda)
\int\limits_{\partial\Omega}\hat{v}(\lambda)\frac{\partial\varphi_{nm}}{\partial\nu}d\sigma+\xi_{nm}}
{\lambda+\lambda^2_{nm}\hat{K}(\lambda)}.
\end{equation}

\textbf{Lemma 1.} If in the problem (\ref{4})---(\ref{6}) $\hat{K}(\lambda)$ has at least one root $\lambda_0\neq 0$
in the domain of holomorphism (we require that this domain exists) then
controllability to rest is impossible; that is, there exists the initial condition $\xi\in H^2(\Omega_0)\cap H^1_0(\Omega_0)$
such that for any $T>0$ and for every control $v\in C(0,T;H^{\frac{3}{2}}(\partial\Omega_0))$
the corresponding solution does not have compact support (in t).

\textbf{Proof.} Let us use the polar coordinates in the integral of the equality (\ref{9}).
Then the equality (\ref{9})
takes the form:
\begin{equation}
\label{10}
\hat{\theta}_{nm}(\lambda)=
\frac{-\mu_n^m\hat{K}(\lambda)
\int\limits_{0}^{2\pi}\hat{v}_0(\lambda,\alpha)e^{im\alpha} d\alpha+\xi_{nm}}
{\sqrt{\pi}R\left(\lambda+\lambda^2_{nm}\hat{K}(\lambda)\right)},
\end{equation}
where $\hat{v}_0(\lambda,\alpha):=\hat{v}(\lambda,R\cos\alpha,R\sin\alpha)$.

The system of functions $\{e^{im\alpha}\}_{m\in Z}$ is an orthogonal basis in $L_2(0,2\pi)$. Thus we can expand
$$
\hat{v}_0(\lambda,\alpha)=\sum\limits_{j=-\infty}^{+\infty}\hat{v}_{0,j}(\lambda)e^{ij\alpha},
$$
where
$$
\hat{v}_{0,j}(\lambda)=\frac{1}{2\pi}\int\limits_0^{2\pi}e^{-ij\alpha}\hat{v}_0(\lambda,\alpha)d\alpha.
$$
Hence we obtain
\begin{equation}
\label{11}
\hat{\theta}_{nm}(\lambda)=
\frac{-\mu_n^m\hat{K}(\lambda)
2\pi\hat{v}_{0,-m}(\lambda) + \xi_{nm}}
{\lambda+\lambda^2_{nm}\hat{K}(\lambda)}.
\end{equation}

We note that if the system is controllable to rest then $\theta_{nm}(t)$, $v_{0,-m}(t)$ have
compact support. Thus $\hat{\theta}_{nm}(\lambda)$ and $\hat{v}_{0,-m}(\lambda)$ are in $PW_+$.

As it follows from the definition of $PW_+$, $\hat{\theta}_{nm}(\lambda)$ is an entire function then
it can not have singularities at the roots of the denominator $\lambda+\lambda^2_{nm}\hat{K}(\lambda)$.
Thus the control function $\hat{v}_{0,-m}(\lambda)$ has to satisfy the following equalities:
\begin{equation}
\label{12}
\hat{v}_{0,-m}(\lambda)=-\frac{1}{2\pi}\frac{\lambda^2_{nm}\xi_{nm}}{\mu^m_n\lambda}
\end{equation}
when $\lambda\neq 0$ is a root of the equation $\lambda+\lambda^2_{nm}\hat{K}(\lambda)=0$.
As $\lambda^2_{nm}=\left(\mu^m_n\right)^2/R^2$ then (\ref{12}) can be rewritten as follows
\begin{equation}
\label{12.1}
\hat{v}_{0,-m}(\lambda)=-\frac{1}{2\pi}\frac{\mu^m_n\xi_{nm}}{R^2\lambda}.
\end{equation}

We note that equalities (\ref{12.1}) can be presented in the following form:
$$
\int\limits_0^T v_{0,-m}(t)e^{-\lambda t}dt=-\frac{1}{2\pi}\frac{\mu^m_n\xi_{nm}}{R^2\lambda}.
$$
The latter equalities are the so-called moment problem.

We record now index $m$. Let $m$, for example, be equal to $1$ and $n$ changes from $1$ to $+\infty$. We get
the subsystem of equalities for values of the function $\hat{v}_{0,-1}(\lambda)$ at such points $\lambda$
that $\lambda+\lambda^2_{n1}\hat{K}(\lambda)=0$:
\begin{equation}
\label{13}
\hat{v}_{0,-1}(\lambda)=-\frac{1}{2\pi}\frac{\mu^1_n\xi_{n1}}{R^2\lambda}.,\quad n=1,2,...\:.
\end{equation}
We note that $\hat{K}(\lambda)$ has a root $\lambda_0\neq 0$ (if $K(t)$ is a series of decreasing exponentials then
$\hat{K}(\lambda)$ has a countable number of roots, see \cite{Ivanov_Sheronova}). Then
applying methods used in \cite{Ivanov} (in which Rouche's theorem was used) we can prove that
there exists a sequence $\{\lambda_n\neq 0\}$ of zeros of
$$
\lambda+\frac{(\mu^1_n)^2}{R^2}\hat{K}(\lambda)
$$
and it is important that this sequence is convergent
to the non-zero complex number.
Let us choose $\xi_{2j+1,1}=0$. Hence $\hat{v}_{0,-1}(\lambda_{2n+1})=0$. As the sequence of zeros
is convergent and $\hat{v}_{0,1}(\lambda)$ is an entire function, then $\hat{v}_{0,-1}(\lambda)\equiv 0$.
We get that for any $n$ all $\xi_{2n,1}$ have to be zero.
But we can always take some of them as non-zero numbers.
Thus we come to the conclusion that
there exists such initial condition $\xi$ that for any control function $v$ controllability to rest
is impossible. Lemma is proved.

We consider the problem (\ref{1})---(\ref{3}). The following theorem is the main result of this article.

\textbf{Theorem.} If the control function $u\in C(0,T;L_2(\Omega))$ in the equation (\ref{1}) has compact support
(in $x,y$) on $\Omega$ and $\hat{K}(\lambda)$ has at least one root $\lambda_0\neq 0$
in the domain of holomorphism
then controllability to rest is impossible; that is, there exists the initial condition
$\xi\in H^2(\Omega)\cap H^1_0(\Omega)$
such that for any $T>0$ and for every control $u\in C(0,T;L_2(\Omega))$
(with compact support on $\Omega$)
the corresponding solution does not have compact support (in t).

\textbf{Proof.} Let $u(t,x,y)$ at each time $t>0$
be an arbitrary function in $L_2(\Omega)$ with compact support $D$.
As support of u is a closed and bounded set then we can consider
a circle $\Omega_0$ which is properly contained in $\Omega$ and does not intersect $D$ (see figure 1).
At each time $t>0$ the solution $\theta(t)\in H^2(\Omega)\cap H^1_0(\Omega)$. We restrict the solution $\theta(t)$
to $\Omega_0$ (in this case $\theta(t)\in H^2(\Omega_0)$) and
consider a new initial boundary-value
problem on $\Omega_0$. In this problem the boundary
condition is equal to the restriction (in the sense of the trace theorem)
of the solution on the boundary of $\Omega_0$. As $\theta(t)\in H^2(\Omega_0)$
this restriction can be computed and at the each time $t>0$ is an element
of $H^{\frac{3}{2}}(\partial\Omega_0)$. We can consider
the restriction of the solution to $\partial\Omega_0$ as
the control. Thus we obtain the boundary control problem on the circle $\Omega_0$.
The solution of this problem exists automatically.
Using the previous lemma we have proved that if controllability to rest of this new problem is
impossible then it is impossible, to stop oscillations of the origin problem (\ref{1})---(\ref{3}).
The theorem is proved.

\section{Some Generalization and Related Topics}
\label{sec:3}

To obtain analogous result for Gurtin-Pipkin equation in the case where the dimension of
$\Omega$ is greater than 2, it is
necessary to use the orthonormalized system of eigenvectors of $A$ defined in $\Omega_0$,
where $\Omega_0$ is a ball. These eigenvectors, for example in $\mathbb{R}^3$, are constructed by means of
spherical harmonics
$$
Y_m^l(\alpha,\varphi),\: 0\leq\alpha\leq\pi,\:
0\leq\varphi< 2\pi,\:
m=0,1,...,\: l=0,\pm1,\pm2,...,\pm m,
$$
while we use functions $e^{im\alpha}$, $0\leq\alpha<2\pi$, $m=0,1,...,$ in case of two-dimensional domains.
Using orthogonal property of spherical harmonics $Y_m^l$ on the unit sphere
we expand $\hat{v}_0(\lambda)$ and after that
all steps of the proof remain unchanged (see proof of the lemma 1). We choose $\mathbb{R}^2$ in present paper
as the clearest and evident way to demonstrate the idea of the proof.

Most likely, it can be proved that if
$$
K(t)=\sum\limits_{j=1}^{N}c_je^{-\gamma_jt}
$$
and the control function does not have compact support properly contained in $\Omega$,
then the problem (\ref{1})---(\ref{3}) is controllable to rest.
We note that if $K(t)=C>0$ then the equation (\ref{1}) is a classical wave equation and
$\hat{K}(\lambda)$ does not have a null, and it is a well-known fact that
using a control contained in a subdomain the system is controllable to rest.
If $K(t)=qe^{-\gamma t}$, $q,\gamma>0$ then the equation (\ref{1}) can be reduced to the equation
$$
\theta_{tt}(t,x,y)-q\Delta\theta(t,x,y)+\gamma\theta_t(t,x,y)=P(t,x,y),
$$
where
$$
P(t,x,y)=\frac{du(t,x,y)}{dt}+\gamma u(t,x,y)
$$
and $P$ can be considered as a new
control. The latter equation is a damped wave equation.
Let us consider the one-dimensional case then, instead of the Laplace operator $\Delta$,
we write the second derivative $\frac{d^2}{dx^2}$:
\begin{equation}
\label{17}
\theta_{tt}(t,x)-q\theta_{xx}(t,x)+\gamma\theta_t(t,x)=P(t,x).
\end{equation}
It is proved (will be published later) that oscillations of the string governed
by the equation (\ref{17}) (let $P(t,x)\equiv 0$) can be stopped
if we apply the control to the end of the string, the second end being fixed.
Apparently, using this fact it can be proved (but it is not proved),
that by means of a control $P(t,x)$ contained in a subsegment (in $x$),
the system is also controllable to rest.

Finally, we also note a link between stability and controllability for the one-dimensional case. If we consider equation
$$
\theta_{tt}-\alpha\theta_{xx}-q\theta_{xx}\ast e^{-\gamma t}=0,
$$
where $\ast$ is convolution and $\alpha>0$, then solutions of this equation are stable if the parameter
$q\in[0,\alpha\gamma]$, and unstable if $q<0$ or $q>\alpha\gamma$. Furthermore, if $q=0$ or $q=\alpha\gamma$
then the system is controllable to rest.

\section{Conclusions}

In this article, we have proved that the system governed by the two-dimensional Gurtin-Pipkin equation
is uncontrollable to rest if the distributed control is supported on the
subdomain which is properly contained in arbitrary bounded domain with a smooth boundary.
In this case, the memory kernel is a twice continuously differentiable function, such that its
Laplace transformation has at least one non-zero root.

\begin{center}
\begin{figure}
  \includegraphics[scale=0.6]{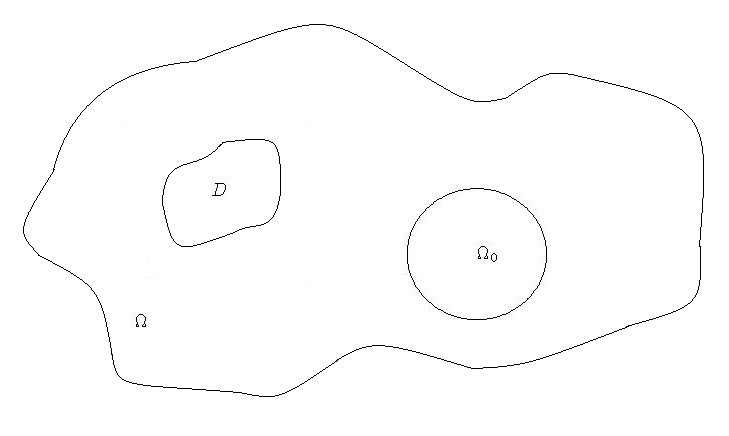}
  \caption{}
\end{figure}
\end{center}

\end{document}